\documentclass{amsart}[12pt]
\usepackage{xypic}
\usepackage{eucal}

\usepackage{amsmath}
\usepackage{amstext}
\usepackage{amssymb}
\usepackage{amsthm}
\usepackage{amscd}

\usepackage{lscape}
\usepackage{longtable}

\usepackage{epsfig}
\input txdtools
\let\et=\etexdraw
\def\etexdraw{\drawbb\et}

\theoremstyle{plain}%default
\newtheorem{thm}{Theorem}[section]
\newtheorem{thm*}{Theorem}
\newtheorem{lem}[thm]{Lemma}
\newtheorem{prop}[thm]{Proposition}
\newtheorem{prop*}[thm*]{Proposition}
\newtheorem{cor}[thm]{Corollary}

\theoremstyle{definition}
\newtheorem{defi}[thm]{Definition}

\newtheorem{ex}[thm]{Example}

\newtheorem{ntn}[thm]{Notation}

\theoremstyle{remark}

\newtheorem{rmk}[thm]{Remark}

\DeclareMathOperator{\Rad}{rad}
\DeclareMathOperator{\height}{ht}
\DeclareMathOperator{\Hom}{Hom}

\DeclareMathOperator{\Ext}{Ext}
\DeclareMathOperator{\Tor}{Tor}
\DeclareMathOperator{\Ann}{ann}
\DeclareMathOperator{\Span}{Span}
\DeclareMathOperator{\grAnn}{gr-ann}

\begin{document}

\def\Ker {\operatorname{Ker}\nolimits}
\def\Im {\operatorname{Im}\nolimits}
\def\Image {\operatorname{Image}\nolimits}
\def\Syz {\operatorname{Syz}\nolimits}
\def\initial {\operatorname{in}\nolimits}
\def\Gin{\operatorname{Gin}\nolimits}
\def\Spec{\operatorname{Spec}\nolimits}
\def\D{\operatorname{D}\nolimits}
\def\V{\operatorname{V}\nolimits}
\def\H {\operatorname{H}\nolimits}
\def\E {\operatorname{E}\nolimits}
\def\V {\operatorname{V}\nolimits}
\def\nE {\operatorname{e}\nolimits}
\def\nV {\operatorname{v}\nolimits}
\def\C {\operatorname{\cal C}\nolimits}
\def\row {\operatorname{row}\nolimits}
\def\column {\operatorname{column}\nolimits}

\title
[$F$-stable submodules] {$F$-stable submodules of top local cohomology modules of Gorenstein rings}
\author{Mordechai Katzman}
%\author{Rodney Y.~Sharp}
\address{Department of Pure Mathematics,
University of Sheffield, Hicks Building, Sheffield S3 7RH, United Kingdom\\
{\it Fax number}: +44-(0)114-222-3769}
\email{M.Katzman@sheffield.ac.uk}

%\subjclass{Primary 13P99, 13D40, 05C69, 05C38}

%\date{\today}

%\keywords{Graded commutative Noetherian ring, graded local cohomology module, infinite set of associated primes.}

%\begin{abstract}
%\end{abstract}

\maketitle
%{\bf\huge Preliminary version-- not for circulation}
\begin{abstract}
This paper applies G.~Lyubeznik's notion of $F$-finite modules to
describe in a very down-to-earth manner certain annihilator submodules of
some top local cohomology modules over Gorenstein rings.
As a consequence we obtain an explicit description of the test ideal of
Gorenstein rings in terms of ideals in a regular ring.
\end{abstract}
%%%%%%%%%%%%%%%%%%%%%%%%%%%%%%%%%%%%%%%%%%%%%%%%%%%%%%%%%%%%%%%%%%%%%%%%%%%%%%%%%%%%%%%%%%%%%
\section{Introduction}\label{Section: Introduction}
%%%%%%%%%%%%%%%%%%%%%%%%%%%%%%%%%%%%%%%%%%%%%%%%%%%%%%%%%%%%%%%%%%%%%%%%%%%%%%%%%%%%%%%%%%%%%

\emph{
Throughout this paper
$(R,\frak{m})$ will denote a regular local ring of characteristic $p$, and $A$ will be a surjective image of $R$.
We also denote the injective hull of $R/m$ with $E$ and for any $R$-module $N$ we write $\Hom_R (N,E)$ as
$N^\vee$.}
We shall always
denote with $f:R\rightarrow R$ the Frobenius map, for which
$f(r)= r^p$ for all $r \in R$ and we shall denote the $e$th iterated Frobenius functor over $R$ with $F^e_R(-)$.
As $R$ is regular, $F^e_R(-)$ is exact (cf.~Theorem 2.1 in \cite{Kunz}.)

For any commutative ring $S$ of characteristic $p$, the skew polynomial ring $S[T; f]$ associated to $S$ and
the Frobenius map $f$
is a non-commutative ring which as a
left $R$-module is freely generated by $(T^i)_{i \geq 0}$,
 and so consists
 of all polynomials $\sum_{i = 0}^n s_i T^i$, where  $n\geq 0$
 and  $s_0,\ldots,s_n \in S$; however, its multiplication is subject to the
 rule
 $$
  Ts = f(s)T = s^pT \quad \mbox{~for all~} s \in S\/.
 $$

Any $A[T; f]$-module  $M$ is a $R[T; f]$-module in a natural way and, as $R$-modules,
$F^e_R(M)\cong RT^e \otimes_R M$.

It has been known for a long time that the local cohomology module $\H^{\dim A}_{\mathfrak{m}A} (A)$
has the structure of an $A[T; f]$-module and this fact has been employed by many authors to study problems
related to tight closure and to Frobenius closure. Recently R.~Y.~Sharp has described in \cite{Sharp}
the parameter test ideal of $F$-injective rings in terms of certain $A[T; f]$-submodules of $\H^{\dim A}_{\mathfrak{m}A} (A)$ and
it is mainly this work which inspired us to look further into the structure of these $A[T; f]$-modules.

The main aim of this paper is to produce a description of the
$A[T; f]$-submodules of $\H^{\dim A}_{\mathfrak{m}A} (A)$ in terms of ideals of $R$ with certain properties.
We first do this when $A$ is a complete intersection.
The $F$-injective case is described by Theorem \ref{Theorem: bijection}
and as a corollary we obtain a description of the parameter test ideal of $A$.
Notice that for Gorenstein rings the test ideal the parameter test ideal coincide
(cf.~Proposition 8.23(d) in \cite{Hochster-Huneke-0} and Proposition 4.4(ii) in \cite{Smith2}.)
We then proceed to describe  the parameter test ideal
in the non-$F$-injective case (Theorem \ref{Theorem: non-$F$-injective}.)
We generalise these results to Gorenstein rings in section \ref{Section: The Gorenstein case}.

%%%%%%%%%%%%%%%%%%%%%%%%%%%%%%%%%%%%%%%%%%%%%%%%%%%%%%%%%%%%%%%%%%%%%%%%%%%%%%%%%%%%%%%%%%%%%
\section{Preliminaries: $F$-finite modules}\label{Section: Preliminaries: $F$-finite modules}
%%%%%%%%%%%%%%%%%%%%%%%%%%%%%%%%%%%%%%%%%%%%%%%%%%%%%%%%%%%%%%%%%%%%%%%%%%%%%%%%%%%%%%%%%%%%%

The main tool used in this paper is
the notion of $F$-modules, and in particular $F$-finite modules.
These were introduced in G.~Lyubeznik's seminal work \cite{Lyubeznik} and  provide
a very fruitful point of view of local cohomology modules in prime characteristic $p$.

One of the tools introduced in \cite{Lyubeznik} is a functor
$\mathcal{H}_{R,A}$ from the category of $A[T; f]$-modules which are Artinian as $A$-modules
to the category of $F$-finite modules.
For any $A[T; f]$-module $M$ which is Artinian as an $A$-module
the $F$-finite structure of $\mathcal{H}_{R,A}(M)$ is obtained as follows.
Let $\gamma : R T \otimes_R M \rightarrow M$ be the $R$-linear map defined by
$\gamma (r T \otimes m) = r T m$; apply the functor $^\vee$ to obtain
$\gamma^\vee : M^\vee \rightarrow F_R(M)^\vee$. Using the isomorphism between
$F_R(M)^\vee$ and $F_R(M^\vee)$ (Lemma 4.1 in \cite{Lyubeznik}) we obtain a map
$\beta: M^\vee \rightarrow F_R(M^\vee)$ which we adopt as a generating morphism of $\mathcal{H}_{R,A}(M)$.

\emph{
We shall henceforth assume that the kernel of the surjection $R\rightarrow A$
is minimally generated by $\mathbf{u}=(u_1, \dots, u_n)$.}
We shall also assume until section \ref{Section: The Gorenstein case} that
$A$ is a complete intersection.
We shall write $u=u_1\cdot \ldots \cdot u_n$ and
for all $t\geq 1$ we let $\mathbf{u}^t R$ be the ideal $u_1^t R +  \dots +  u_n^t R$.

To obtain the results in this paper we shall need to understand the $F$-finite module structure of
$$\mathcal{H}_{R,A}\left( \H^{\dim A}_{\mathfrak{m}A} (A) \right) \cong \H^{\dim R - \dim A}_{\mathbf{u}R}(R) ;$$
this has generating root
$$\frac{R}{\mathbf{u} R} \xrightarrow[]{u^{p-1}} \frac{R}{\mathbf{u}^p R} $$
(cf.~Remark 2.4 in \cite{Lyubeznik}.)

\begin{defi}
Define $\mathcal{I}(R,\mathbf{u})$ to be the set of all ideals $I\subseteq R$ containing $(u_1, \dots , u_n)R$
with the property that
$$u^{p-1} \left(I + \mathbf{u}R\right) \subseteq I^{[p]}+\mathbf{u}^pR .$$
\end{defi}

\begin{lem}\label{Lemma: generating morphisms}
Consider the $F_R$-finite $F$-module
$M=\H^n_{\mathbf{u} R}(R)$
with generating root
$$\frac{R}{\mathbf{u} R} \xrightarrow[]{u^{p-1}} \frac{R}{\mathbf{u}^pR} .$$

\begin{enumerate}
\item[(a)]
For any $I\in \mathcal{I}(R,\mathbf{u})$ the  $F_R$-finite module
with generating root
$$\frac{I+ \mathbf{u}R}{\mathbf{u} R} \xrightarrow[]{u^{p-1}} \frac{I^{[p]}+\mathbf{u}^p R}{\mathbf{u}^p R} \cong F_R\left( \frac{I+\mathbf{u}R}{\mathbf{u}R} \right) $$
is an $F$-submodule of  $M$ and every $F_R$-finite $F$-submodule of  $M$ arises in this way.
\item[(b)]
For any $I\in \mathcal{I}(R,\mathbf{u})$ the  $F_R$-finite module
with generating morphism
$$\frac{R}{I+\mathbf{u}R} \xrightarrow[]{u^{p-1}} \frac{R}{I^{[p]}+\mathbf{u}^p R} \cong F_R\left( \frac{R}{I+\mathbf{u} R}\right) $$
is an $F$-module quotient of  $M$ and every $F_R$-finite $F$-module quotient  of  $M$ arises in this way.
\end{enumerate}
\end{lem}

\begin{proof}
\noindent (a) For any $I\in \mathcal{I}(R,\mathbf{u})$, the map
$$\frac{I+\mathbf{u}R}{\mathbf{u} R} \xrightarrow[]{u^{p-1}} \frac{I^{[p]}+\mathbf{u}^p R}{\mathbf{u}^p R} $$
is well defined and is injective; now the first statement follows from
Proposition 2.5(a) in \cite{Lyubeznik}.
If $N$ is any $F_R$-finite $F$-submodule of  $M$, the root of $N$ is a submodule of the root of $M$, i.e.,
the root of $N$ has the form $(I+\mathbf{u}R)/\mathbf{u}R$ for some ideal $I\subseteq R$ (cf.~\cite{Lyubeznik}, Proposition 2.5(b))
and the structure morphism of $N$ is induced by that of $M$, i.e., by multiplication by $u^{p-1}$, so we must have
$u^{p-1} I \subseteq I^{[p]}+\mathbf{u}^p R $, i.e., $I\in \mathcal{I}(R,\mathbf{u})$.

\noindent (b)
For any $I\in \mathcal{I}(R,\mathbf{u})$, the map
$$\frac{R}{I+\mathbf{u}R} \xrightarrow[]{u^{p-1}} \frac{R}{I^{[p]}+\mathbf{u}^p R} \cong F_R\left( \frac{R}{I+\mathbf{u}R}\right) $$
is well defined and we have the following  commutative diagram with exact rows
\begin{equation}\label{CD3}
\xymatrix{
0 \ar@{>}[r]^{} &
\displaystyle\frac{I+\mathbf{u}R}{\mathbf{u}R} \ar@{>}[r]^{} \ar@{>}[d]^{u^{p-1}} &
\displaystyle\frac{R}{\mathbf{u}R} \ar@{->}[r]^{} \ar@{>}[d]^{u^{p-1}} &
\displaystyle\frac{R}{I+\mathbf{u}R} \ar@{>}[r]^{} \ar@{>}[d]^{u^{p-1}}  &
0 \\
0 \ar@{>}[r]^{} &
\displaystyle F_R\left(\frac{I+\mathbf{u}R}{\mathbf{u}R}\right) \ar@{>}[r]^{} \ar@{>}[d]^{u^{p(p-1)}}                  &
\displaystyle  F_R\left(\frac{R}{\mathbf{u}R}\right) \ar@{>}[r]^{} \ar@{-}[d]^{u^{p(p-1)}} &
\displaystyle F_R\left(\frac{R}{I+\mathbf{u}R}\right) \ar@{>}[r]^{} \ar@{>}[d]^{u^{p(p-1)}}  &
0 \\
& \vdots &\vdots&\vdots&\\
}
\end{equation}
Taking direct limits of the vertical maps we obtain an exact sequence
$0 \rightarrow M^{\prime} \rightarrow M \rightarrow M^{\prime\prime} \rightarrow 0$
which establishes the first statement of (b).

Conversely, if $M^{\prime\prime}$ is a $F$-module quotient of $M$, say, $M^{\prime\prime}\cong M/M^{\prime}$
for some $F$-submodule $M^{\prime}$ of  $M$ use (a) to find a generating root of $M^{\prime}$ of the form
$$\frac{I+\mathbf{u}R}{\mathbf{u} R} \xrightarrow[]{u^{p-1}} \frac{I^{[p]}+\mathbf{u}^pR}{\mathbf{u}^pR} $$
for some $I\in \mathcal{I}(R,\mathbf{u})$. Looking again at the direct limits of the vertical maps in (\ref{CD3})
we establish the second statement of (b).

\end{proof}

\begin{defi}
For all $I\in \mathcal{I}(R,\mathbf{u})$ we define $\mathcal{N}(I)$ to be the $F$-module quotient of
$\H^n_{\mathbf{u}R}(R)$ with generating morphism
$$\frac{R}{I+\mathbf{u}R} \xrightarrow[]{u^{p-1}} \frac{R}{I^{[p]}+\mathbf{u}^pR} \cong F_R\left( \frac{R}{I+\mathbf{u}R}\right) .$$
\end{defi}

\begin{lem}\label{Lemma 2}
Assume that $R$ is complete.
Let $H$ be an Artinian $A[T; f]$-module and write $M=\mathcal{H}_{R,A}(H)$.
Let $N$ be a homomorphic image of $M$ with generating morphism $N_0$.
Then $N_0^\vee$ is an $A[T; f]$-submodule of $H$ and  $N\cong \mathcal{H}_{R,A}(N_0^\vee)$.
\end{lem}
\begin{proof}
Notice that $M$ (and hence $N$) are $F$-finite modules (cf.~\cite{Lyubeznik}, Theorems 2.8 and 4.2).
Let $N_0$  be root of $N$ and $M_0$ a root of $M$
so that we have a commutative diagram with exact rows
\begin{equation*}\label{CD1}
\xymatrix{
M_0 \ar@{->>}[r]^{} \ar@{->}[d]^{\mu} & N_0 \ar@{->}[r]^{}  \ar@{->}[d]^{\nu} & 0 \\
F_R(M_0) \ar@{->>}[r]^{}              & F_R(N_0) \ar@{->}[r]^{}   & 0
}
\end{equation*}
where the vertical arrows are generating morphisms.
Apply the functor $\Hom(-,E)$ to the commutative diagram above
to obtain the following commutative diagram with exact rows
\begin{equation*}\label{CD2}
\xymatrix{
0 \ar@{->}[r]^{}  & F_R(N_0)^\vee \ar@{->}[r]^{} \ar@{->}[d]^{\nu^\vee} &  F_R(M_0)^\vee  \ar@{->}[d]^{\mu^\vee}\\
0 \ar@{->}[r]^{}  & N_0^\vee \ar@{->}[r]^{}  & M_0^\vee
}
\end{equation*}
and recall that  $M_0$ is isomorphic to $H^\vee$ (cf.~\cite{Lyubeznik}, Theorem 4.2).
Since $R$ is complete, $\left(H^\vee\right)^\vee\cong H$ and we immediately see that
$N_0^\vee$ is a $R$-submodule of $H$.  We now show that $N_0^\vee$ is an $A[T;f]$ submodule of $H$ by showing that
$T N_0^\vee \subseteq N_0^\vee$.

The construction of the functor $\mathcal{H}_{R,A}(-)$ is such that
for any $h\in H \cong  M_0^\vee$, $T h$ is the image of $T \otimes_R h$ under the map
$$F_R(M_0)^\vee  \xrightarrow[]{\mu^\vee} M_0^\vee $$
and so for $h\in N_0^\vee$,
$T h$ is the image of $T \otimes_R h$ under the map
$$F_R(N_0)^\vee  \xrightarrow[]{\nu^\vee} N_0^\vee $$
and hence $T  h\in  N_0^\vee$.

Now the fact that $N\cong \mathcal{H}_{R,A} \left(N_0^\vee\right)$ follows the construction of the
functor $\mathcal{H}_{R,A}(-)$.
\end{proof}

\begin{ntn}
Let $M$ be a left $A[T, f]$-module.
We shall write $AT^\alpha M$ for the $A$-module generated
by $T^\alpha M$. Note that $AT^\alpha M$ is a left $A[T, f]$-module.
We shall also write $\displaystyle M^\bigstar=\bigcap_{\alpha\geq 0} AT^\alpha M$.
\end{ntn}

\begin{lem}\label{Lemma 4}
Assume that $R$ is complete.
Let $H$ be an $A[T;f]$-module and assume that $H$ is $T$-torsion-free.
Let $I,J \subseteq A$ be ideals. If, for some $\alpha\geq 0$,
$$ A T^\alpha \Ann_H {I A[T; f]}  = A T^\alpha \Ann_H {J A[T; f]}  $$
then $\Ann_H { I A[T; f]} =\Ann_H { J A[T; f]} $.
\end{lem}
\begin{proof}
Both $A T^\alpha \Ann_H {I A[T; f]} $ and $A T^\alpha \Ann_H {J A[T; f]}  $
are left $A[T; f]$-submodules.
Now for every $T$-torsion-free $A[T;f]$-module $M$, and every ideal $K\subseteq A$, if
$$ \left(\bigoplus_{i\geq 0} K T^i\right) AT^\alpha M = \left(\bigoplus_{i\geq 0} K T^{i+\alpha}\right) M$$
vanishes then so does
$$\left(\bigoplus_{i\geq 0} K^{[p^\alpha]} T^{i+\alpha}\right) M=
\left(\bigoplus_{i\geq 0} T^{\alpha} K T^{i}\right) M=
T^{\alpha} \left(\bigoplus_{i\geq 0}  K T^{i}\right) M$$
and since $M$ is $T$-torsion-free,
$$ \left(\bigoplus_{i\geq 0}  K T^{i}\right) M=0 .$$
We deduce that $\grAnn A T^\alpha M=\grAnn M$.
Now
%Lemma 3.2 in \cite{Sharp} implies that
$$\grAnn A T^\alpha \big(\Ann_H {I A[T; f]} \big) = \grAnn \Ann_H {I A[T; f]}  ,$$
$$\grAnn A T^\alpha \big(\Ann_H {J A[T; f]}  \big)= \grAnn \Ann_H {J A[T; f]} $$
and Lemma 1.7 in  \cite{Sharp} shows that $\Ann_H{ I A[T; f]} =\Ann_H { J A[T; f]} $.
\end{proof}

%%%%%%%%%%%%%%%%%%%%%%%%%%%%%%%%%%%%%%%%%%%%%%%%%%%%%%%%%%%%%%%%%%%%%%%%%%%%%%%%%%%%%%%%%%%%%%%%%%%%%%%%%
\section{The $A[T;f]$ module structure of top local cohomology modules of $F$-injective Gorenstein rings}
%%%%%%%%%%%%%%%%%%%%%%%%%%%%%%%%%%%%%%%%%%%%%%%%%%%%%%%%%%%%%%%%%%%%%%%%%%%%%%%%%%%%%%%%%%%%%%%%%%%%%%%%%

\begin{defi}
As in \cite{Smith2} we say that an ideal $I\subseteq A$ is an \emph{$F$-ideal} if
$\Ann_{\H^{\dim(A)}_{\frak{m} A}(A)} I$ is a left $A[T; f]$-module, i.e., if
$\Ann_{\H^{\dim(A)}_{\frak{m} A}(A)} I=\Ann_{\H^{\dim(A)}_{\frak{m} A}(A)} IA[T; f]$.
\end{defi}

\begin{thm}\label{Theorem 1}
Assume that $R$ is complete.
Consider the $F_R$-finite $F$-module
$ M=\H^n_{\mathbf{u}R}(R)$
with generating root
$$\frac{R}{\mathbf{u}R} \xrightarrow[]{u^{p-1}} \frac{R}{\mathbf{u}^pR} $$
and consider the Artinian $A[T; f]$ module $H=\H^{\dim(A)}_{\frak{m} A}(A)$.
Let $N$ be a homomorphic image of $M$.

\begin{enumerate}
\item[(a)] $M=\mathcal{H}_{R,A}(-)(H)$ and has generating root
$H^\vee \cong R/\mathbf{u}R \xrightarrow[]{u^{p-1}}  R/\mathbf{u}^p R \cong F_R(H^\vee)$.
\item[(b)]
If $N$ has generating morphism
$$\frac{R}{I+\mathbf{u}R} \xrightarrow[]{u^{p-1}} \frac{R}{I^{[p]}+\mathbf{u}^pR} $$
then $I A$ is an $F$-ideal, $N \cong \mathcal{H}_{R,A}\big( \Ann_H {IA[T;f]} \big)$.
If, in addition, $H$ is $T$-torsion free then
$\grAnn \Ann_H {IA[T;f]}=IA[T;f]$ and $I$ is radical.

\item[(c)] Assume that $H$ is $T$-torsion free (i.e., $H_r=H$ in the terminology of \cite{Lyubeznik}).
For any ideal $J\subset R$, the $F$-finite module $\mathcal{H}_{R,A}\big(\Ann_H {JA[T;f]} \big)$
has generating morphism
$$\frac{R}{I+\mathbf{u}R} \xrightarrow[]{u^{p-1}} \frac{R}{I^{[p]}+\mathbf{u}^pR} $$
for some ideal $I\in \mathcal{I}(R,\mathbf{u})$ with $\Ann_H { I A[T; f]} =\Ann_H { J A[T; f]} $.
\end{enumerate}
\end{thm}

\begin{proof}
The first statement is
a restatement of the discussion at the beginning of section \ref{Section: Preliminaries: $F$-finite modules}.
%an immediate consequence of Example 4.8 in \cite{Lyubeznik}: indeed $H^\vee \cong \Ext^n_R(R/I,R) \cong R/\mathbf{u}R$ and  $ F_R(H^\vee) \cong \Ext^n_R(R/I^{[p]},R) \cong R/\mathbf{u}^p R $ and the map $\Ext^n_R(R/I,R) \rightarrow \Ext^n_R(R/I^{[p]},R)$ induced by the quotient map $R/I^{[p]} \rightarrow R/I$ arises from the map of Koszul complexes $\mathcal{K}_\bullet(\mathbf{u}; R) \rightarrow \mathcal{K}_\bullet(\mathbf{u}^p; R)$. This map induces a map of cohomologies $\Ext^n_R(R/I,R) \rightarrow \Ext^n_R(R/I^{[p]},R)$ which is given by multiplication by $u$. %\marginpar{\small Show that the generating morphism is given by multiplication by $u^{p-1}$!}

Notice that Lemma \ref{Lemma: generating morphisms} implies that $N$ must have a generating morphism
of the form given in (b) for some $I\in \mathcal{I}(R,\mathbf{u})$.

Since $A$ is Gorenstein, $H$ is an injective hull of $A/\frak{m}A$ which we denote $\overline{E}$.
Lemma \ref{Lemma 2} implies that
$ N \cong \mathcal{H}_{R,A}\left( L \right)$ where
$\displaystyle L= \left( \frac{R}{I+\mathbf{u}R} \right)^\vee$ is a $A[T; f]$-submodule of $H=\overline{E}$.
But
\begin{eqnarray*}
\left( \frac{R}{I+\mathbf{u}R} \right)^\vee&=& \Ann_E \left(I+\mathbf{u}R\right) \\
&=& \Ann_{\left( \Ann_{\mathbf{u}R} E \right)} I  \\
&=& \Ann_{\overline{E}} I   .
\end{eqnarray*}
But  $L$ is a $A[T; f]$-submodule of $\overline{E}$ and so $I A$ is an $F$-ideal
%implies  that
%$L\subseteq \Ann_{I A[T; f]} \overline{E}$.
%Since $I \Ann_{I A[T; f]} \overline{E}=0$,  $L\supseteq \Ann_{I A[T; f]} \overline{E}$ and
%we obtain
and $L= \Ann_{\overline{E}} {I A[T; f]}$.
Also,
\begin{eqnarray*}
\left( 0 :_R \Ann_{\overline{E}} {I A[T; f]} \right) &= & \left( 0 :_R \Ann_E I \right) \\
&=&\left( 0 :_R (R/I)^\vee \right) \\
&=&\left( 0 :_R (R/I) \right) \\
&=& I
\end{eqnarray*}
(where the third equality follows from 10.2.2 in \cite{Brodmann-Sharp})
If $H$ is $T$-torsion free,
Proposition 1.11 in \cite{Sharp} implies that $I=\grAnn \Ann_{\overline{E}} {I A[T; f]}$
and Lemma 1.9 in \cite{Sharp} implies that $I$ is radical.

To prove part (c) we recall Lemma \ref{Lemma: generating morphisms} which states
that $\mathcal{H}_{R,A}\big(\Ann_H {JA[T;f]} \big)$
has generating morphism
$$\frac{R}{I+\mathbf{u}R} \xrightarrow[]{u^{p-1}} \frac{R}{I^{[p]}+\mathbf{u}^pR} $$
for some $I\in \mathcal{I}(R,\mathbf{u})$ and we need only  show that
$\Ann_H { I A[T; f]} =\Ann_H { J A[T; f]} $.

Part (b) implies that
$\mathcal{H}_{R,A}\big(\Ann_H {JA[T;f]} \big) = \mathcal{H}_{R,A} \big(\Ann_H {IA[T;f]} \big)$
for some $I\in \mathcal{I}(R,\mathbf{u})$
and
Theorem 4.2 (iv) in \cite{Lyubeznik} implies
$$ \bigcap_{i=0}^\infty A T^i \big( \Ann_H {JA[T;f]}  \big)= \bigcap_{i=0}^\infty A T^i \big( \Ann_H {IA[T;f]}  \big) $$
and since $H$ is Artinian there exists an $\alpha\geq 0$ for which
$A T^\alpha \big( \Ann_H {JA[T;f]} \big)=  A T^\alpha \big( \Ann_H {IA[T;f]}  \big) $
and the result follows from Lemma \ref{Lemma 4}.
\end{proof}

\begin{rmk}
Theorem \ref{Theorem 1}  can provide an easy way to show
that $H=\H^{\dim(A)}_{\frak{m} A}(A)$ is \emph{not} $T$-torsion free. As an example consider
$R=\mathbb{K}[\![ x,y, a,b ]\!]$, $u=x^2a-y^2b$ and $A=R/uR$.
Its easy to verify that $(x,y,a^2)R\in \mathcal{I}(R,x^2a-y^2b)$
when $\mathbb{K}$ has characteristic 2,
and we deduce that $\H^{3}_{(x,y,a,b) A}(A)$ is not  $T$-torsion free.
\end{rmk}

\begin{thm}\label{Theorem2}
Assume that $R$ is complete and that $\H^{\dim(A)}_{\frak{m} A}(A)$ is $T$-torsion free.
\begin{enumerate}
\item[(a)] For all $A[T;f]$-submodules $L$ of $\H^{\dim(A)}_{\frak{m} A}(A)$,
$$L^\bigstar=\bigcap_{i=0}^\infty A T^i L$$
has the form
%\marginpar{Define $L^*$ as in Lyubeznik}
$A T^{\alpha} M$ where $\alpha\geq 0$ and $M$ is a special annihilator submodule in the terminology of \cite{Sharp}.
\item[(b)] The set $\big\{  \mathcal{N}(I) \,|\, I\in \mathcal{I}(R,\mathbf{u}) \big\}$ is finite.
\end{enumerate}
\end{thm}
\begin{proof}
\noindent (a)
Let $L$ be a  $A[T;f]$-submodule of $\H^{\dim(A)}_{\frak{m} A}(A)$.
Pick a $I\in \mathcal{I}(R,\mathbf{u})$
such that $\mathcal{N}(I)=\mathcal{H}_{R,A}\left(  L \right)$.
Now use part (b) of Theorem \ref{Theorem 1} and deduce that
$\mathcal{N}(I)\cong \mathcal{H}_{R,A}\left( \Ann_H {I A[T;f]}   \right)$.
Now the result follows from Theorem 4.2 (iv) in \cite{Lyubeznik}.

\noindent (b)
Theorem \ref{Theorem 1}(b) implies that
$$\big\{  \mathcal{N}(I) \,|\, I\in \mathcal{I}(R,\mathbf{u}) \big\}=
\left\{  \mathcal{H}_{R,A} \left(\Ann_{\H^{\dim A}_{\frak{m} A} (A) } {I A[T; f]} \right)\,|\, I\in \mathcal{I}(R,\mathbf{u}) \right\} ;$$
now Corollary 3.11 and Proposition 1.11 in \cite{Sharp} imply that the set on the right is finite.

\end{proof}

The following Theorem reduces the problem of classifying all
$F$-ideals of $A$ (in the terminology of \cite{Smith2}) or all
special $\H^{\dim(A)}_{\frak{m} A}(A)$-ideals  (in the terminology of \cite{Sharp}) in the case where
$A$ is an $F$-injective complete intersection, to problem of
determining the set $\mathcal{I}(R,\mathbf{u})$.

\begin{thm}\label{Theorem: bijection}
Assume $H:=\H^{\dim(A)}_{\frak{m} A}(A)$ is $T$-torsion free and let $\mathcal{B}$ be the set of all
$H$-special $A$-ideals (cf. \S 0 in \cite{Sharp})
\begin{enumerate}
\item[(a)] The map $\Psi: \mathcal{I}(R,\mathbf{u}) \rightarrow \mathcal{B}$ given by
$\Psi(I)=I A$ is a bijection.
\item[(b)] There exists a unique minimal element $\tau$ in
$\left\{I \,|\, I\in \mathcal{I}(R,\mathbf{u})  , \ \height I A >0  \right \}$
and that $\tau$ is a parameter-test-ideal for $A$.
\item[(c)] $A$ is $F$-rational if and only if  $\mathcal{I}(R,\mathbf{u})=\{ 0, R \}$.
\end{enumerate}
\end{thm}
\begin{proof}
\noindent (a)
Assume first that $R$ is complete.
Theorem \ref{Theorem 1}(b) implies that $\Psi$ is well defined, i.e.,
$\Psi(I)\in \mathbf{B}$ for all $I\in\mathcal{I}(R,\mathbf{u})$, and, clearly, $\Psi$ is injective.
The surjectivity of $\Psi$ is a consequence of Theorem \ref{Theorem 1}(c).

Assume now that $R$ is not complete, denote completions with $\widehat{\phantom{h}}$ and write
$\widehat{H}=\H^{\dim(\widehat{A})}_{\frak{m} \widehat{A}}(\widehat{A})$.
If $I$ is a $\widehat{H}$-special $\widehat{A}$-ideal, i.e.,
if there exists an $\widehat{A}[T; f]$-submodule $N\subseteq \widehat{H}$ such that
$\grAnn N=I \widehat{A}[T; f]$ then
$I=(0 :_{\widehat{A}} N)$ (cf.~Definition 1.10 in \cite{Sharp}).
But recall that $\widehat{H}=H$ and $N$ is a $A[T; f]$-submodule of $H$; now
$I=(0 :_{\widehat{A}} N)=(0 :_A N) \widehat{A} $.
If we let  $\widehat{\mathcal{B}}$ be the set of
$\H^{\dim(\widehat{A})}_{\frak{m} \widehat{A}}(\widehat{A})$-special
$\widehat{A}$-ideals, we have a bijection $ \Upsilon : \mathcal{B} \rightarrow \widehat{\mathcal{B}}$
mapping $I$ to $I\widehat{A}$.
This also shows that all ideals in $\mathcal{I}(\widehat{R},\mathbf{u})$  are expanded
from $R$, and now since $\widehat{R}$ is faithfully flat over $R$,
we deduce that all ideals in $\mathcal{I}(\widehat{R},\mathbf{u})$ have the form
$I \widehat{R}$ for some $I \in \mathcal{I}(R,\mathbf{u})$.
We now obtain a chain of bijections
$$ \mathcal{I}(R,\mathbf{u}) \longleftrightarrow
\mathcal{I}(\widehat{R},\mathbf{u}) \longleftrightarrow
\widehat{\mathcal{B}} \longleftrightarrow
\mathcal{B} .$$

\noindent (b)
This is immediate from (a) and Corollary 4.7 in \cite{Sharp}.

\noindent (c)
If $A$ is $F$-rational,  $\H^{\dim(A)}_{\frak{m} A}(A)$ is a simple $A[T;f]$-module
(cf.~Theorem  2.6 in \cite{Smith3}) and the only $H$-special $A$-ideals must be $0$ and $A$.
The bijection established in (a) implies now $\mathcal{I}(R,\mathbf{u})=\{ 0, R \}$.

Conversely, if $\mathcal{I}(R,\mathbf{u})=\{ 0, R \}$, part (b) of the Theorem implies that $1\in A$ is a parameter-test-ideal,
i.e., for all systems of parameters $\mathbf{x}=(x_1, \dots, x_d)$ of $A$,
$\left( \mathbf{x} A \right)^*=\left( \mathbf{x} A \right)^F= \mathbf{x} A $
where the second equality follows from the fact that $\H^{\dim(A)}_{\frak{m} A}(A)$ is $T$-torsion free.

\end{proof}

%%%%%%%%%%%%%%%%%%%%%%%%%%%%%%%%%%%%%%%%%%%%%%%%%%%%%%%%%%%%%%%%%%%%%%%%%%%%%%%%%%%%%%%%%%
\section{Examples}
%%%%%%%%%%%%%%%%%%%%%%%%%%%%%%%%%%%%%%%%%%%%%%%%%%%%%%%%%%%%%%%%%%%%%%%%%%%%%%%%%%%%%%%%%%
%%%%% Try x^4+ y*z*t^2+ y*z^2*t+y^2*z*t: F-injective(p=5), not F-rational. Only special prime is (x,y,z,t).

Throughout this section $\mathbb{K}$ will denote a field of prime characteristic.

\begin{ex}\label{Example 1}
Let $R$ be the localization of $\mathbb{K}[ x,y ]$ at $(x,y)$, $u=xy$ and $A=R/uR$.
Then $\H^1_{xy R}(R)=\mathcal{H}_{R,A}\big(\H^1_{xA+yA}(A)\big)$ ought to have four proper $F$-finite $F$-submodules
corresponding to the elements $0$, $xR$, $yR$ and $xR+yR$ of $\mathcal{I}(R,xy)$.

We verify this by giving an explicit
description the $A[T;f]$-module structure of
$$H:=\H^1_{xA+yA}(A)\cong
\lim_{\longrightarrow}
\left(
\frac{A}{(x-y) A } \xrightarrow[]{x-y}
\frac{A}{(x-y)^2 A} \xrightarrow[]{x-y}
\frac{A}{(x-y)^3 A} \xrightarrow[]{x-y} \dots
\right)
$$
First notice that in $H$, for all $n\geq 1$ and $0<\alpha\leq n$,
$x^\alpha + (x-y)^n A= x +  (x-y)^{n-\alpha+1}$ and
$y^\alpha + (x-y)^n A= y +  (x-y)^{n-\alpha+1}$ so
$H$ is the $\mathbb{K}$-span of $\left\{x + (x-y) A\right\} \cup X \cup Y \cup U$ where
$$
\begin{array}{l}
X=\big\{ x + (x-y)^n A  \,|\, n\geq 2 \big\},\\
Y=\big\{ y + (x-y)^n A  \,|\, n\geq 2 \big\},\\
U=\big\{ 1 + (x-y)^n A  \,|\, n\geq 1 \big\}
\end{array}
$$
and notice also that the action of the Frobenius map $f$ on $H$ is such that
$ T\left( x^\alpha + (x-y)^n A\right) =  x^{\alpha p} + (x-y)^{n p} A $
and
$ T\left( y^\alpha + (x-y)^n A \right) =  y^{\alpha p} + (x-y)^{n p} A $ for all $\alpha\geq 0$.

Next notice that any $A[T, f]$-submodule $M$ of $H$ which contains an element
$1 + (x-y)^n A \in U$
must coincide with $H$: for $1\leq m<n$ we have
$(x-y)^{n-m}\left(  1 + (x-y)^n A \right) =  (x-y)^{n-m} + (x-y)^n A = 1 + (x-y)^m A$,
whereas for $m>n$, pick
an $e\geq 0$ such that $n p^e>m$, write
$$
T^e \left( 1 + (x-y)^n A \right)=  1 +  (x-y)^{n p^e}A \in M$$
and  use the previous case ($m<n$) to deduce that $1 + (x-y)^m A\in M$.
Since now $U\subseteq M$, we see that $M=H$.

We now show that there are only three non-trivial $A[T, f]$-submodules of $H$, namely $\Span_\mathbb{K} X$
and $\Span_\mathbb{K} Y$, and $\Span_\mathbb{K} \left\{x + (x-y) A\right\} \cup X$.
By symmetry, it is enough to show that,
if $M$ is an $A[T, f]$-submodule of $H$ and
$x + (x-y)^n A\in M$ for some $n\geq 2$,
then $X\subset M$.
If $1\leq m<n$,
$$x^{n-m} \left( x + (x-y)^n A \right)=
x^{n-m+1} + (x-y)^n A=
x + (x-y)^{n-(n-m)}A=
x + (x-y)^m A $$
whereas, if $m>n\geq 2$, pick an $e\geq 0$ such that $n p^e-p^e+1>m$ and write
$$
T^e \left( x + (x-y)^n A \right)=
x^{p^e} + (x-y)^{n p^e} A =
x + (x-y)^{n p^e-p^e+1} A \in M$$
and  using the previous case ($m<n$) we deduce that
$x + (x-y)^m A \in M$.
\end{ex}

\begin{ex}
%%%%% Try x^2*y+x*y*z+z^3: has special primes (x,y,z) and (z,x), jacobian=(z,x), F-injective, not F-rational.
Let $R$ be the localization of $\mathbb{K}[x,y,z]$ at ${\frak m}=(x,y,z)$, $u=x^2y+xyz+z^3$ and $A=R/uR$.
Fedder's criterion (cf.~Propositon 2.1 in \cite{Fedder}) implies that $A$ is $F$-pure,
and Lemma 3.3 in \cite{Fedder} implies that
the $A[T;f]$ module $\H^1_{{\frak m}A}(A)$ is $T$-torsion-free.

Here $\mathcal{I}(R,u)$ contains
the ideals $0$, $xR+zR$ and $xR+yR+zR$.
We deduce that $A$ is not $F$-rational and that its parameter-test-ideal is $xR+zR$.
Also, Theorem \ref{Theorem: bijection}(b) implies that the only proper ideals in $\mathcal{I}(R,u)$
% of positive height
are the ones listed above.
\end{ex}

\begin{ex}
%%%%% Try x^2*y+x*y*z+z^3: has special primes (x,y,z) and (z,x), jacobian=(z,x), F-injective, not F-rational.
Let $R$ be the localization of $\mathbb{K}[ x,y,z ]$ at ${\frak m}=(x,y,z)$
and assume that $\mathbb{K}$ has characteristic $2$. Let
$u=x^3+y^3+z^3+xyz$ and $A=R/uR$.
Notice that we can factor $u=(x+y+z)(x^2+y^2+z^2+xy+xz+yz)$.
Fedder's criterion implies that $A$ is $F$-pure,
and Lemma 3.3 in \cite{Fedder} implies that
the $A[T;f]$ module $\H^1_{{\frak m}A}(A)$ is $T$-torsion-free.

Here
\begin{eqnarray*}
\mathcal{I}(R,u) & \supseteq & \big\{ 0, (x+y+z)R, (x^2+y^2+z^2+xy+xz+yz)R, \\
&& (x+z,y+z)R, (x+y+z,y^2+yz+z^2)R, \\
&& (x,y,z)R \big\} .
\end{eqnarray*}
The images in $A$ of the first three ideals have height zero while
the images in $A$ of the fourth and fifth ideals have height $1$.
Using \ref{Theorem: bijection}(b)
we conclude that that the parameter test-ideal of $A$ is a sub-ideal of
$$J=(x+z,y+z)A \cap (x+y+z,y^2+yz+z^2)A=(x^2+yx, y^2+xz, z^2+xy)A .$$
But this ideal defines the singular locus of $A$ and Theorem 6.2 in \cite{Hochster-Huneke-1}
implies that the parameter test-element of $A$ contains $J$, so $J$ is the parameter test-ideal of $A$.
\end{ex}

\section{The non-$F$-injective case}

In this section we extend the results of the previous section to the case where $A$ is not $F$-injective.
First we produce a criterion for the $F$-injectivity of $A$.

\begin{defi}
Define
$$\mathcal{I}_0(R,\mathbf{u})=\left\{ L\in \mathcal{I}(R,\mathbf{u}) \,|\,
u^{(p-1)(1+p+\dots+p^{e-1})}  \in L^{[p^e]} + \mathbf{u}^{p^e} R \text{ for some } e\geq 1 \right\} .$$
\end{defi}

\begin{prop}
\begin{enumerate}
\item[(a)] For any  $L\in \mathcal{I}(R,\mathbf{u})$, $\mathcal{N}(L)=0$ if and only if $L\in \mathcal{I}_0(R,\mathbf{u})$.
\item[(b)] $\H^{\dim(A)}_{\frak{m} A}(A)$ is $T$-torsion free if and only if $\mathcal{I}_0(R,\mathbf{u})=\{R\}$.
\end{enumerate}
\end{prop}
\begin{proof}
\noindent(a)
Recall that the $F$-finite module $\mathcal{N}(L)$ has generating morphism
$$\frac{R}{L+\mathbf{u}R} \xrightarrow[]{u^{p-1}}
\frac{R}{L^{[p]}+\mathbf{u}^pR} \cong F_R\left( \frac{R}{L+\mathbf{u}R}\right) .$$
Proposition 2.3 in \cite{Lyubeznik} implies that
$\mathcal{N}(L)=0$ if and only if for some $e\geq 1$
the composition
$$\frac{R}{L+\mathbf{u}R}
\xrightarrow[]{u^{p-1}}
\frac{R}{L^{[p]}+\mathbf{u}^pR}
\xrightarrow[]{u^{(p-1)p}}
\frac{R}{L^{[p]}+\mathbf{u}^{p^2}R}
\dots
\xrightarrow[]{u^{(p-1)p^{e-1}}}
\frac{R}{L^{[p]}+\mathbf{u}^{p^e}R}$$
vanishes, i.e., if and only if
$u^{(p-1)(1+p+\dots+p^{e-1})}  \in L^{[p^e]} + \mathbf{u}^{p^e} R $
for some $e\geq 1$.

\noindent(b)
Write $H=\H^{\dim(A)}_{\frak{m} A}(A)$.
If $H$ is $T$-torsion free, the existence of the bijection
described in Theorem \ref{Theorem: bijection}(a) implies that for any non-unit
$L\in\mathcal{I}_0(R,\mathbf{u})$, $\Ann_H {L A[T;f]}\neq \Ann_H {A[T;f]}=0$.
Theorem \ref{Theorem 1}(b) implies $\mathcal{N}(L)\cong \mathcal{H}_{R,A}\big( \Ann_H {L A[T;f]} \big)$
so $\mathcal{H}_{R,A}\big( \Ann_H {L A[T;f]} \big)=0$. But Theorem 4.2(ii) in \cite{Lyubeznik} now implies
that $\Ann_H {L A[T;f]}$ is nilpotent, a contradiction.

Assume now that $H$ is not $T$-torsion free, i.e., $H_n\neq 0$.
The short exact sequence
$$0 \rightarrow H_n \rightarrow H \rightarrow H/H_n \rightarrow 0 $$
yields the short exact sequence
$$0 \rightarrow \big(H/H_n\big)^\vee \rightarrow \frac{R}{\mathbf{u} R} \rightarrow H_n^\vee \rightarrow 0 .$$
Notice that as the functor $\Hom(-, E)$ is faithful, $H_n^\vee \neq 0$, and so
$H_n^\vee\cong R/I$ for some ideal $\mathbf{u} R \subseteq I\varsubsetneq R$.
Now $\mathcal{H}_{R,A}\big( H_n \big)$
is the $F$-finite quotient of $H$ with generating morphism
$$\frac{R}{I} \xrightarrow[]{u^{p-1}} \frac{R}{I^{[p]}} $$
and this vanishes because of Theorem 4.2(ii) in \cite{Lyubeznik}, i.e., $I\in \mathcal{I}_0(R,\mathbf{u})$.
\end{proof}

We now describe the parameter test ideal of $A$.
Henceforth we shall always denote $\H^{\dim(A)}_{\frak{m} A}(A)$  with $H$.

\begin{thm}\label{Theorem: non-$F$-injective}
Assume that $R$ is complete.
The parameter test ideal of $A$ is given by
$$\bigcap\big\{ I\in \mathcal{I}(R,\mathbf{u}) \,|\, \height IA>0 \big\}. $$
\end{thm}

\begin{proof}
Write $\overline{\tau}$ for the parameter test ideal of $A$ and let $\tau$ be its pre-image in $R$.
Recall that $\overline{\tau}$ is an $F$-ideal (Proposition 4.5 in \cite{Smith2},)
i.e., $\Ann_H \overline{\tau}$ is an $A[T;f]$-submodule of $H$, and
$\mathcal{H}_{R,A}\big( \Ann_H \overline{\tau} \big)$ has generating morphism
$$\left( \Ann_H \overline{\tau} \right)^\vee \xrightarrow[]{u^{p-1}} F_R\left(\left( \Ann_H \overline{\tau} \right)^\vee\right) .$$
But
$$\left( \Ann_H \overline{\tau} \right)^\vee \cong
\left( \left( A/\overline{\tau} \right)^\vee \right)^\vee \cong
R/(\tau+\mathbf{u}R) $$
so the generating morphism of $\mathcal{H}_{R,A}\big( \Ann_H \overline{\tau} \big)$
is
$$R/(\tau+\mathbf{u}R) \xrightarrow[]{u^{p-1}} R/(\tau^{[p]}+\mathbf{u}^p R) $$
and so we must have $\tau\in \mathcal{I}(R,\mathbf{u})$.

As $A$ is Cohen-Macaulay,
$\overline{\tau}=\displaystyle (0 :_A 0^*_H) $ (cf.~Proposition 4.4 in \cite{Smith2}.)

By Theorem \ref{Theorem 1}(b), for each $I\in \mathcal{I}(R,\mathbf{u})$, the ideal $IA$ is an $F$-ideal and,
if $\height I>0$,
$\displaystyle \Ann_H IA= \Ann_H I A[T;f] \subseteq 0^*_H $
and so
$$\displaystyle \overline{\tau}=(0 :_A 0^*_H) \subseteq
\bigcap\big\{ (0 :_A \Ann_H IA) \,|\, IA\in \mathcal{I}(R,\mathbf{u}), \height IA>0 \big\} .$$
But $H$ is an injective hull of $A/\mathfrak{m}A$ so
$$ ( 0 :_A \Ann_H IA) = \left( 0 :_A \Hom(A/IA, H) \right)= (0 :_A A/IA)=IA $$
and
$$\overline{\tau}\subseteq \bigcap\big\{IA \,|\, IA\in \mathcal{I}(R,\mathbf{u}), \height IA>0 \big\}  .$$
But as $\overline{\tau}$ is one of the ideals in this intersection, we obtain
$\displaystyle \overline{\tau}=
\bigcap\big\{ IA\in \mathcal{I}(R,\mathbf{u}) \,|\, \height IA>0 \big\}$.
\end{proof}

%%%%%%%%%%%%%%%%%%%%%%%%%%%%%%%%%%%%%%%%%%%%%%%%%%%%%%%%%%%%%%%%%
\section{The Gorenstein case}\label{Section: The Gorenstein case}
%%%%%%%%%%%%%%%%%%%%%%%%%%%%%%%%%%%%%%%%%%%%%%%%%%%%%%%%%%%%%%%%%

In this section we generalise the results so far to the case where $A$ is Gorenstein.

Write $\delta=\dim R - \dim A$ and $\overline{E}=E_A(A/\mathfrak{m}A)$.
Local duality implies
$\displaystyle
\Ext^{\delta}_R(A, R)=
\H^{\dim A}_\mathfrak{m}(A)^\vee \cong
\Hom\left(\H^{\dim A}_{\mathfrak{m}A}(A),\overline{E} \right)$
and since $A$ is Gorenstein this is just $A=R/\mathbf{u}R$.

Now
$\Ext^{\delta}_R \left( R/\mathbf{u}R, A \right)\cong R/\mathbf{u} R$,
$\Ext^{\delta}_R \left( R/\mathbf{u}^pR, A \right)\cong R/\mathbf{u}^p R$
and
$\displaystyle \mathcal{H}_{R,A}\big( \H^{\dim A}_{\mathfrak{m} A} \big)= \H^\delta_{\mathfrak{m}}(R)$
has generating morphism
$R/\mathbf{u} \rightarrow R/\mathbf{u}^p R$
given by multiplication by some element of $R$ which we denote
$\varepsilon(\mathbf{u})$ (this is unique up to multiplication by a unit.)
Unlike the complete intersection case, the map
$R/\mathbf{u} \xrightarrow[]{\varepsilon(\mathbf{u})}  R/\mathbf{u}^p R$
may not be injective, i.e., this generating morphism of  $\H^\delta_{\mathfrak{m}}(R)$
is not a \emph{root}.
However, if define
$$\displaystyle K_\mathbf{u}:=
\bigcup_{e\geq 0} \big( \mathbf{u}^{p^{e+1}} R :_R \varepsilon(\mathbf{u})^{1+p+\dots+p^e} \big)$$
we obtain a root
$R/K_\mathbf{u} \xrightarrow[]{\varepsilon(\mathbf{u})}  R/K_\mathbf{u}^{[p]}$
(cf.~Proposition 2.3 in \cite{Lyubeznik}.)

We now extend naturally our definition of $\mathcal{I}(R, \mathbf{u})$ when $A$ is Gorenstein as follows.
\begin{defi}
If $A=R/\mathbf{u}R$ is Gorenstein we define
$\mathcal{I}(R, \mathbf{u})$ to be the set of all ideals $I$ of $R$ containing $K_\mathbf{u}$
for which
$\varepsilon(\mathbf{u}) I \subseteq I^{[p]}$.
\end{defi}

Now a routine modification of the proofs of the previous sections gives the following two theorems.

\begin{thm}\label{Theorem: the Gorenstein case I}
Assume $A$ is Gorenstein and that $H^{\dim A}_{\mathfrak{m} A}(A)$ is $T$-torsion-free.
\begin{enumerate}
\item[(a)] The map $I \mapsto IA$ is a bijection between $\mathcal{I}(R, \mathbf{u})$
and the $A$-special $H^{\dim A}_{\mathfrak{m} A}(A)$-ideals.
\item[(b)] There exists a unique minimal element $\tau$ in
$\left\{I \,|\, I\in \mathcal{I}(R,\mathbf{u})  , \ \height I A >0  \right \}$
and that $\tau$ is a parameter-test-ideal for $A$.
\item[(c)] $A$ is $F$-rational if and only if  $\mathcal{I}(R,\mathbf{u})=\{ 0, R \}$.
\end{enumerate}
\end{thm}

\begin{thm}\label{Theorem: the Gorenstein case II}
Assume that $R$ is complete and that $A$ is Gorenstein.
The parameter test ideal of $A$ is given by
$$\bigcap\big\{ I\in \mathcal{I}(R,\mathbf{u}) \,|\, \height IA>0 \big\}. $$
\end{thm}

%Finally we make a connection between these results and Fedder's criterion
%(Propositon 2.1 in \cite{Fedder}) which states that $A$ is $F$-injective if and only if
%$u^{p-1} \mathfrak{m} \notin  \mathfrak{m}^{[p]}$.

\section*{Acknowledgement}
My interest in the study of local cohomology modules as modules over
skew polynomial rings was aroused by many interesting conversations with Rodney Sharp,
during one of which I learnt about Example  \ref{Example 1}.

\end{document}